\documentclass{article}

\usepackage{amsmath}
\usepackage{amssymb}
\usepackage{amsthm}
\usepackage{amscd}

\usepackage{overcite}

\makeatletter
\renewcommand\thetable{\@Roman\c@table}
\makeatother

\newcommand{\Z}{{\mathbb Z}}
\newcommand{\ket}[2]{#1^{\vert #2 \rangle}}
\newcommand{\bra}[2]{#1^{\langle #2 \vert}}
\newtheorem{theorem}{Theorem}[]
\newtheorem{lemma}[theorem]{Lemma}

\newtheorem{proposition}[theorem]{Proposition}
\newtheorem{corollary}[theorem]{Corollary}

\theoremstyle{definition}
\newtheorem{definition}[theorem]{Definition}
\newtheorem{example}[theorem]{Example}

\makeatletter
\def\eqnarray{%
  \stepcounter{equation}%
  \let\@currentlabel=\theequation
  \global\@eqnswtrue
  \global\@eqcnt\z@
  \tabskip\@centering
  \let\\=\@eqncr
  $$\halign to \displaywidth\bgroup\@eqnsel\hskip\@centering
  $\displaystyle\tabskip\z@{##}$&\global\@eqcnt\@ne
  \hfil$\displaystyle{{}##{}}$\hfil
  &\global\@eqcnt\tw@$\displaystyle\tabskip\z@{##}$\hfil
  \tabskip\@centering&\llap{##}\tabskip\z@\cr}
\makeatother

\def\geh{\mathfrak{g}}

\def\cd{\cdots}
\def\ot{\otimes}

\def\veps{\varepsilon}
\def\vphi{\varphi}

\def\Z{{\mathbb Z}}
\def\et{\tilde{e}}
\def\ft{\tilde{f}}
\def\qed{\hfill$\square$\par}

\title{Factorization of Combinatorial $R$ matrices
and Associated Cellular Automata}

\author{
Goro Hatayama\thanks{
Institute of Physics, University of Tokyo, Komaba, Tokyo 153-8902, Japan},
Atsuo Kuniba,$\hspace{-1.2mm}^*$
and Taichiro Takagi\thanks{
Department of Applied Physics, National Defense Academy,
Yokosuka 239-8686, Japan}
}

\date{}

\begin{document}

\maketitle
\begin{center}
{\it Dedicated to Professor Rodney Baxter on the occasion
of his sixtieth birthday }
\end{center}

\begin{abstract}
Solvable vertex models in statistical mechanics give rise to 
soliton cellular automata at $q=0$ in a 
ferromagnetic regime.
By means of the crystal base theory we study a class of such 
automata associated with non-exceptional quantum affine algebras 
$U'_q(\hat{\geh}_n)$.
Let $B_l$ be the crystal of the $U'_q(\hat{\geh}_n)$-module corresponding to 
the $l$-fold symmetric fusion of the vector representation.
For any crystal of the form $B = B_{l_1} \ot \cd \ot B_{l_N}$,
we prove that the combinatorial $R$ matrix
$B_M \ot B \xrightarrow{\sim} B \ot B_M$ is factorized into a
product of Weyl group operators in a certain domain if
$M$ is sufficiently large.
It implies the factorization of certain transfer matrix at $q=0$,
hence the time evolution in the associated cellular automata.
The result generalizes the ball-moving algorithm 
in the box-ball systems.\\[.5eM]
\end{abstract}

\section{Introduction}\label{sec:intro}

The box-ball systems \cite{TS,T,TTM,TNS} are important 
examples of soliton cellular automata.
They are discrete dynamical systems whose time evolution is 
expressed as a certain motion of balls along the one dimensional array of boxes.
Their integrability has been understood by
the ultradiscretization \cite{TTMS} of  
classical integrable systems (soliton equations).
In the recent works \cite{HKT}${}^,$\cite{FOY,HHIKTT}
it was revealed that the box-ball systems may also be viewed as  
quantum integrable systems at $q=0$.
Here by quantum integrable systems we mean the ones whose 
integrability is guaranteed by  the Yang-Baxter equation \cite{B}, and 
$q$ is the deformation parameter in the relevant quantum group.
In fact, the box-ball systems are identified with 
a $q \rightarrow 0$ limit of some two-dimensional solvable vertex models,
where the role of time evolution is played by 
the action of a row transfer matrix.
The simplest example is the original Takahashi-Satsuma automaton \cite{TS},
whose classical origin is the discrete Lotka-Volterra equation \cite{TTMS}
and the quantum origin is the fusion six-vertex model.
Here is an example of the automaton time evolution.

\vskip0.2cm
\noindent
$\hspace*{4cm}\cdots 1112211211111111 \cdots$\\ 
$\hspace*{4cm}\cdots 1111122121111111 \cdots$\\
$\hspace*{4cm}\cdots 1111111212211111 \cdots$\\
$\hspace*{4cm}\cdots 1111111121122111 \cdots$
\vskip0.2cm

\noindent
One regards $1$ as an empty box and $2$ as a box containing a ball.
At each time step one moves every ball once 
starting from the leftmost one. 
The rule is that the ball goes to the nearest right empty box.
One easily finds that the sequence of $ l$ balls propagate 
stably to the right with velocity $l$ unless it interacts with other balls.
By regarding  such patterns as (ultradiscrete) solitons, 
the above figure illustrates
how the larger soliton overtakes the 
smaller one with a phase shift.

In terms of the fusion six vertex model at $q=0$, 
the above figure corresponds to the configuration: 

{\small 
\setlength{\unitlength}{0.9mm}
\begin{picture}(130,55)(-2,7)

\multiput(5,10)(12,0){10}{\line(0,1){8}}
\multiput(2,14)(12,0){10}{\line(1,0){6}}

\multiput(5,25)(12,0){10}{\line(0,1){8}}
\multiput(2,29)(12,0){10}{\line(1,0){6}}

\multiput(5,40)(12,0){10}{\line(0,1){8}}
\multiput(2,44)(12,0){10}{\line(1,0){6}}

\put(-4.5,43 ){111}
\put(-4.5,28){111}
\put(-4.5,13){111}

\put(8.5,43 ){112}
\put(8.5,28){111}
\put(8.5,13){111}

\put(20.5,43 ){122}
\put(20.5,28){111}
\put(20.5,13){111}

\put(32.5,43 ){112}
\put(32.5,28){112}
\put(32.5,13){111}

\put(44.5,43 ){111}
\put(44.5,28){122}
\put(44.5,13){111}

\put(56.5,43 ){112}
\put(56.5,28){112}
\put(56.5,13){112}

\put(68.5,43 ){111}
\put(68.5,28){122}
\put(68.5,13){111}

\put(80.5,43 ){111}
\put(80.5,28){112}
\put(80.5,13){112}

\put(92.5,43 ){111}
\put(92.5,28){111}
\put(92.5,13){122}

\put(104.5,43 ){111}
\put(104.5,28){111}
\put(104.5,13){112}

\put(116.5,43 ){111}
\put(116.5,28){111}
\put(116.5,13){111}

\put(4.2,50){2}
\put(4.2,35){1}
\put(4.2,20){1}
\put(4.2,5){1}

\put(16.2,50){2}
\put(16.2,35){1}
\put(16.2,20){1}
\put(16.2,5){1}

\put(28.2,50){1}
\put(28.2,35){2}
\put(28.2,20){1}
\put(28.2,5){1}

\put(40.2,50){1}
\put(40.2,35){2}
\put(40.2,20){1}
\put(40.2,5){1}

\put(52.2,50){2}
\put(52.2,35){1}
\put(52.2,20){2}
\put(52.2,5){1}

\put(64.2,50){1}
\put(64.2,35){2}
\put(64.2,20){1}
\put(64.2,5){2}

\put(76.2,50){1}
\put(76.2,35){1}
\put(76.2,20){2}
\put(76.2,5){1}

\put(88.2,50){1}
\put(88.2,35){1}
\put(88.2,20){2}
\put(88.2,5){1}

\put(100.2,50){1}
\put(100.2,35){1}
\put(100.2,20){1}
\put(100.2,5){2}

\put(112.2,50){1}
\put(112.2,35){1}
\put(112.2,20){1}
\put(112.2,5){2}

\end{picture}
}

\noindent
The left and right boundaries are to be understood as  $1$ or $111$ everywhere.
This is a configuration of the fusion six vertex model in which 
the quantum space is spin $1/2$ ($1$ or $2$) and the auxiliary space is 
spin $3/2$ ($111, 112, 122$ or $222$).
At $q=0$ only some selected vertex configurations
have non-zero Boltzmann weights and the transfer matrix
yields a deterministic evolution of the spins on one row to another.
The vertex configurations in the above figure 
are the non-zero ones, and form an example of the 
combinatorial ($q=0$) $R$ matrix, which will be a main 
subject in this paper.
Let $T_M$ denote the row transfer matrix at $q=0$ corresponding to 
the spin $M/2$ auxiliary space.
The above example corresponds to $T_3$.
Actually it can be shown that the ball-moving algorithm 
coincides with the action of $T_M$ with sufficiently large $M$.
In the above example $T_M = T_3$ holds for any $M \ge 3$.

The coincidence of  the 
ultradiscrete limit of soliton equations 
and the $q \rightarrow  0$ limit of  vertex models 
is an interesting phenomenon in various respects.
{}From a statistical mechanical point of view, it roughly means that 
in those solvable vertex models,
the profile of low-lying excitations over the ferromagnetic ground state 
at $q=0$ admits an exact description by (ultradiscrete) soliton equations.
{}From a mathematical point of view it leads to 
a systematic generalization \cite{HKT} 
by means of the quantum affine algebras and the crystal base theory
\cite{K1}.
See section 4 of the reference \cite{HKOTY} for several examples of 
the scattering of ultradiscrete solitons.

Now we turn to a general setting, where the six vertex model is 
replaced with a solvable vertex model associated with the quantum affine algebra
$U'_q(\hat{\geh}_n)$.
In this paper we treat the non-exceptional case
$\hat{\geh}_n = A^{(1)}_n, A^{(2)}_{2n-1}, A^{(2)}_{2n},
B^{(1)}_n, C^{(1)}_n, D^{(1)}_n$ and $D^{(2)}_{n+1}$.
The box-ball systems correspond to  \cite{HHIKTT} 
$\hat{\geh}_n = A^{(1)}_n$.
The Boltzmann weights are trigonometric functions satisfying the
Yang-Baxter equation.
The row transfer matrix is specified by the 
auxiliary space $V_M$ and  the quantum space
$\cd \ot V_{l_j}\ot V_{l_{j+1}}\ot \cd$.
Here $V_M$ denotes the $M$-fold symmetric fusion of the
vector representation.
We suppose a ferromagnetic boundary condition, namely,
the spins in the distance $\vert j \vert \gg 1$
are all equal to some prescribed element in $V_{l_j}$.
At $q=0$  the transfer matrix
yields a deterministic evolution of the spins on one row to another.

To analyze such a situation we invoke the crystal base theory \cite{K1}.
Let $B_l$ be the crystal of $V_l$.
It is a finite set listed in Appendix \ref{app:crystal}
endowed with the action of Kashiwara
operators $\et_i, \ft_i : B_l \rightarrow B_l \cup \{0\}$ for
$0 \le i  \le n$.
Let $\delta_l[a] \in B_l$ be the special element as in (\ref{eq:delta}).
The states in the automaton are the elements
$\cd \ot b_j \ot b_{j+1} \ot \cd \in
\cd \ot B_{l_j} \ot B_{l_{j+1}} \ot \cd$
obeying the boundary condition
$b_j = \delta_{l_j}[a_k],\; \vert j \vert \gg 1$.
Here $a_k$ is specified in Table \ref{tab:2} with
(\ref{eq:extendedai}), and $k \in {\mathbb Z}$ is a label of the
boundary condition at our disposal.
Let $T_M$ denote
the $q=0$ transfer matrix, or $M$-th time evolution in the automaton;
$T_M: \cd \ot B_{l_j} \ot B_{l_{j+1}} \ot \cd
\rightarrow \cd \ot B_{l_j} \ot B_{l_{j+1}} \ot \cd$.
The map $T_M( \cdots \ot b_{j} \ot
b_{j+1} \ot \cdots ) = \cdots \ot b'_{j} \ot
b'_{j+1} \ot \cdots$
is induced by the isomorphism of crystals
$R: B_M \ot B \xrightarrow{\sim} B \ot B_M$ with
$B = \cd \ot B_{l_j} \ot B_{l_{j+1}}\ot \cd$
according to  (\ref{eq:evolution}).
We call  $R$ the {\em combinatorial $R$ matrix}.
It is obtained by successive applications of the elementary ones
$B_{M} \ot B_{l_j} \xrightarrow{\sim}
B_{l_j} \ot B_{M}$.
The boundary condition matches the known properties like
$\delta_M[a_k] \ot \delta_{l_j}[a_k] \stackrel{\sim}{\mapsto}
\delta_{l_j}[a_k] \ot \delta_M[a_k]$.
When $M$ gets large, $T_M$ stabilizes to a certain map, which we denote by $T$.
In the box-ball system terminology ($\hat{\geh}_n = A^{(1)}_n$ case),
this corresponds to the
boxes with inhomogeneous capacities $\{ l_j \}$
and the carrier with infinite capacity.

The main result of this note is Theorem \ref{th:Conj9}, which states  that
the isomorphism
$R: B_M \ot B \xrightarrow{\sim} B \ot B_M$ with
$B = B_{l_1} \ot \cd \ot B_{l_N}$ is expressed as
\begin{equation*}
R = (\sigma_B \ot \sigma) P
S_{i_{k+d}}\cd S_{i_{k+2}}S_{i_{k+1}}
\end{equation*}
in the domain  $B_M[a_k]\ot B \subset B_M \ot B$ (\ref{eq:B[]def})
if $M$ is sufficiently large.
Here $S_i$ is the Weyl group operator\cite{K3}
 (\ref{eq:Sdef}) acting on $B_M \ot B$.
$P(u \ot v) = v \ot u$ is the transposition, $\sigma$ and $\sigma_B$ are the operators
corresponding to the Dynkin diagram automorphism  described around
Proposition \ref{pr:sigma}.
The above result on $R$  reveals
the factorization of the time evolution in the automaton:
(Corollary \ref{cor:Th21})
\begin{equation*}
T = \sigma_B
S_{i_{k+d}}\cd S_{i_{k+2}}S_{i_{k+1}},
\end{equation*}
where $S_i$ here is the one acting on
$\cd \ot B_{l_j} \ot B_{l_{j+1}}\ot \cd$.
See the explanation after Corollary \ref{cor:Th21}.
Such a decomposition is by no means evident
{}from the defining relation (\ref{eq:evolution}).
Note that $T$ is a translation in the sense that
the product $\sigma r_{i_{k+d}}\cdots r_{i_{k+1}}$ ($r_i$ is a simple reflection)
is so in the extended affine Weyl group, if
$\sigma$ is interpreted as the Diagram automorphism acting
on the  weight lattice.
According to the factorization,
one can consider a finer time evolution  ${\mathcal T}_m$ (\ref{eq:takagi})
that includes the original one as ${\mathcal T}_{k+td} = T^t\, (t \ge 0)$.
For $\hat{\geh}_n =  A^{(1)}_n$,  the
change {}from ${\mathcal T}_m(p)$ to ${\mathcal T}_{m+1}(p)$ agrees with
the original ball-moving algorithm in the box-ball systems,
where one touches only the balls with a fixed color.
In particular when  $\forall l_j = 1$, our
Definition \ref{def:Def19}
provides a representation theoretical interpretation of the
earlier observation\cite{HIK}.
For $\hat{\geh}_n \neq A^{(1)}_n$,
the automaton corresponding to $\forall l_j = 1$ with $a_k = 1$
has been introduced previously \cite{HKT}.
The formula (\ref{eq:takagi}) in principle
provides a simple algorithm to compute the refined time
evolution for general  $\l_j$ and $a_k$ in an analogous way to the
ball-moving procedure for the $A^{(1)}_n$ case.
The data  $d \in {\mathbb Z}_{\ge 1}$ and $i_k \in I$ are
specified in Table \ref{tab:2}.
Curiously, they have stemmed
{}from the study\cite{KMOU,KMOTU1,KMOTU2} of Demazure crystals \cite{K2}.
It will be interesting to investigate the present result in the light of
the works  \cite{C,KT,KR,KMOTU2,NY,S}.

\section{Crystals}\label{sec:crystals}
Let $\hat{\geh}_n = A^{(1)}_n (n \ge 1), A^{(2)}_{2n-1} (n \ge 3),
A^{(2)}_{2n} (n \ge 2), B^{(1)}_n (n \ge 3), C^{(1)}_n (n \ge 2),
D^{(1)}_n (n \ge 4)$ and $D^{(2)}_{n+1} (n \ge 2)$.
For each $\hat{\geh}_n$ and $l \in {\mathbb Z}_{\ge 1}$,
the $U'_q(\hat{\geh}_n)$ crystal $B_l$
has been constructed \cite{KMN2,KKM} except for $C^{(1)}_n$ with $l$ odd.
As for $C^{(1)}_n$, $B_l$ here is $B^{1,l}$ in the paper \cite{HKKOT}.
The finite set $B_l$ and
the actions of $\et_i,\ft_i: B_l \rightarrow B_l \cup \{0 \}$
for $i \in I= \{0, 1, \ldots, n\}$ (crystal structure) have been defined.
We employ the same notation as \cite{KKM,HKKOT} and
quote the set $B_l$ in Appendix \ref{app:crystal}.
In particular for $\hat{\geh}_n \neq A^{(1)}_n$,
we will always assume the convention  
\[
x_a = \overline{x}_i,\; 
\overline{x}_a = x_i,\; 
\overline{a}=i \quad \text{if} \quad
a = \overline{i},\; 1 \le i \le n.
\]
Let us recall some other notations and the tensor product rule.
\begin{equation*}
\vphi_i(b) = \max\{j \mid \ft_i^jb \neq 0 \},\qquad
\varepsilon_i(b) = \max\{j \mid \et_i^jb \neq 0 \}.
\end{equation*}
For two crystals $B$ and $B'$, the tensor product
$B\ot B'$ is defined as the set
$B\ot B'=\{b_1\ot b_2\mid b_1\in B,b_2\in B'\}$ with
the actions of $\et_{i}$ and $\ft_{i}$ specified by
\begin{eqnarray}
\et_{i}(b_1\ot b_2)&=&\left\{
\begin{array}{ll}
\et_{i}b_1\ot b_2&\mbox{ if }\vphi_i(b_1)\ge\veps_i(b_2)\\
b_1\ot \et_{i}b_2&\mbox{ if }\vphi_i(b_1) < \veps_i(b_2),
\end{array}\right. \label{eq:ot-e}\\
\ft_{i}(b_1\ot b_2)&=&\left\{
\begin{array}{ll}
\ft_{i}b_1\ot b_2&\mbox{ if }\vphi_i(b_1) > \veps_i(b_2)\\
b_1\ot \ft_{i}b_2&\mbox{ if }\vphi_i(b_1)\le\veps_i(b_2).
\end{array}\right. \label{eq:ot-f}
\end{eqnarray}
Here $0\ot b$ and $b\ot0$ are understood to be $0$.
Consequently one has
\begin{equation}
\left.
\begin{array}{ll}
\vphi_i(b_1\ot b_2) &= \vphi_i(b_2) + (\vphi_i(b_1) -\veps_i(b_2))_+,
\\
\veps_i(b_1\ot b_2) &= \veps_i(b_1) + (\veps_i(b_2) -\vphi_i(b_1))_+,
\end{array}
\right.
\label{eq:phi}
\end{equation}
where the symbol $(x)_+$ is defined by $(x)_+ = \max(x,0)$.
The Weyl group operator $S_i \; (i \in I)$ is defined by \cite{K3}
\begin{equation}\label{eq:Sdef}
S_i b = \left\{
\begin{array}{ll}
\ft_i^{\vphi_i(b)-\veps_i(b) }b
&\mbox{ if } \vphi_i(b) \ge \veps_i(b),\\
\et_i^{\veps_i(b) - \vphi_i(b) }b
&\mbox{ if }\vphi_i(b) \le \veps_i(b).
\end{array}\right.
\end{equation}
$S_i$ satisfies the Coxeter relations \cite{K3}.
For two crystals $B$ and $B'$, we let $P: B\ot B' \rightarrow B' \ot B$
denote  the transposition $P(u \ot v) = v \ot u$.
It is easy to check
\begin{equation}\label{eq:Lem25}
S_i P = P S_i \quad \text{ for any } i \in I.
\end{equation}

For two crystals of the form
$B = B_{l_1} \ot \cd \ot B_{l_N}$ and
$B' = B_{l'_1} \ot \cd \ot B_{l'_{N'}}$,
the  tensor products $B \ot B'$ and $B' \ot B$
are isomorphic, i.e., they have the same crystal structure.
The isomorphism $R: B \ot B' \xrightarrow{\sim} B' \ot B$
is called the combinatorial $R$ matrix \cite{KMN1}.
(In this paper we do not consider the energy associated with $R$.)
It is obtained by a successive application of the elementary ones
$B_{l_i} \ot B_{l'_j} \xrightarrow{\sim} B_{l'_j} \ot B_{l_i}$.
We will use the same symbols
$\et_i, \ft_i, \veps_i, \vphi_i, S_i, P$ and $R$
irrespective of the crystals that they act.

Let $\Lambda_i$ denote a fundamental weight and let
$P_{cl} = \bigoplus_{i \in I} {\mathbb Z}\Lambda_i$ be the
classical weight lattice.
(See Section 3.1 of the paper \cite{KMN1} for a precise treatment.)
We define a linear map
$\sigma: P_{cl} \rightarrow P_{cl}$ as in
the rightmost column of Table \ref{tab:1}.
It is a Dynkin diagram automorphism.
When $\hat{\geh}_n \neq C^{(1)}_n$ (resp. $\hat{\geh}_n = C^{(1)}_n$), it
agrees with the one introduced after Corollary 4.6.3
of the paper \cite{KMN1} with  $B = B_l$  (resp. $B = B_{2l}$) for any $l$.
\begin{table}[h]
\caption{The data $\sigma'$ and $\sigma$}\label{tab:1}
\begin{center}
\begin{tabular}{l|c|c}
$\hat{\geh}_n$ & $\sigma'$ on $B_l$ & $\sigma$ on $P_{cl}$ \\ \hline
$A^{(1)}_n$ & $a \rightarrow a-1$ & $\Lambda_{a} \rightarrow \Lambda_{a-1}$ \\
$A^{(2)}_{2n-1}$ & $1 \leftrightarrow \overline{1}$ &
$\Lambda_0 \leftrightarrow \Lambda_1$ \\
$A^{(2)}_{2n}$ & $id$ &
$id$ \\
$B^{(1)}_{n}$ & $1 \leftrightarrow \overline{1}$ &
$\Lambda_0 \leftrightarrow \Lambda_1$ \\
$C^{(1)}_{n}$ & $id$ &
$id$ \\
$D^{(1)}_{n}$ & $1 \leftrightarrow \overline{1},
n \leftrightarrow \overline{n}$ &
$\Lambda_0 \leftrightarrow \Lambda_1, \Lambda_n \leftrightarrow \Lambda_{n-1}$ \\
$D^{(2)}_{n+1}$ & $id$ &
$id$ \\
\end{tabular}
\end{center}
\end{table}
We also let $\sigma$ act on the index set $I$ by the rule
$i' = \sigma(i) \Leftrightarrow \sigma(\Lambda_i) = \Lambda_{i'}$.
For $A^{(1)}_n$, the letters $a, a - 1$ should be interpreted mod $n+1$.

We also introduce the data $d \in {\mathbb Z}_{>0}$, and the
sequences
$i_d, \ldots, i_1 \in I$ and $a_d, \ldots, a_0$ for each algebra
as in Table \ref{tab:2}.
\begin{table}[h]
\caption{The data $d, i_k$ and $a_k$}\label{tab:2}
\begin{center}%
\small{%
\begin{tabular}{l|c|c|c}
$\hat{\geh}_n$ & $d$ & $i_d, \ldots, i_1$ &
$a_d, \ldots, a_0$ \\ \hline
$A^{(1)}_n$ & $n$ & $1,2,\ldots, n$ & $1,2,\ldots, n\!+\!1$ \\
$A^{(2)}_{2n\!-\!1}$ & $2n\!-\!1$
& $n\!-\!1,\ldots,3,2, {0,1 \brace 1,0}, 2,3, \ldots,n$
&$\overline{n},\ldots,\overline{3},\overline{2},{1 \brace \overline{1}},2,\ldots,n\!-\!1,n,\overline{n}$ \\
$A^{(2)}_{2n}$ & $2n$
& $n\!-\!1,\ldots,2,1,0,1,2,\ldots,n$
& $\overline{n},\ldots, \overline{2},\overline{1},1,2,\ldots,n,\overline{n}$ \\
$B^{(1)}_{n}$ & $2n\!-\!1$
& $n\!-\!1,\ldots,3,2,{0,1 \brace 1,0}, 2,3, \ldots,n$
&$\overline{n},\ldots,\overline{3},\overline{2},{1 \brace \overline{1}},2,\ldots,n\!-\!1,n,\overline{n}$ \\
$C^{(1)}_{n}$ & $2n$
& $n\!-\!1,\ldots,2,1,0,1,2,\ldots,n$
& $\overline{n},\ldots, \overline{2},\overline{1},1,2,\ldots,n,\overline{n}$ \\
$D^{(1)}_{n}$ & $2n\!-\!2$
& $n,n\!-\!2,\ldots,2,{0,1 \brace 1,0}, 2, \ldots,n\!-\!2,n$
&$n,\overline{n\!-\!1},\ldots,\overline{2},{1 \brace \overline{1}},2,\ldots,n\!-\!1,\overline{n}$ \\
$D^{(2)}_{n\!+\!1}$ & $2n$
& $n\!-\!1,\ldots,2,1,0,1,2,\ldots,n$
& $\overline{n},\ldots, \overline{2},\overline{1},1,2,\ldots,n,\overline{n}$ \\
\end{tabular}
}
\end{center}
\end{table}
Here $a_k$'s
are taken {}from the letters appearing in the description of the crystals $B_l$
in Appendix \ref{app:crystal}.
In the third and the fourth columns of Table \ref{tab:2},
the symbols ${0,1 \brace 1,0}$ and ${1 \brace \overline{1} }$ mean that
either the simultaneous upper choice or the  simultaneous lower choice are
allowed.

\begin{proposition}\label{pr:sigma}
For any $l\in {\mathbb Z}_{>0}$  the operator
$\sigma' :=  S_{i_1} \cd S_{i_d} : B_l \rightarrow B_l$
is a bijection having  the properties:
\begin{equation*}
\sigma'  \ft_i = \ft_{\sigma(i)}\sigma' , \quad
\sigma'  \et_i = \et_{\sigma(i)}\sigma'.
\end{equation*}
Moreover, the action of $\sigma'$  is  explicitly given by the second column of
Table \ref{tab:1}.
\end{proposition}
The second column of Table \ref{tab:1} specifies the transformation of
those letters labeling the elements of $B_l$.
See Appendix \ref{app:crystal}.
For example in $A^{(1)}_n$ case,
$\sigma'((x_1, \ldots, x_{n+1})) = (x_2,\ldots, x_{n+1},x_1)$
in terms of the coordinates.
Similarly for $A^{(2)}_{2n-1}$,
$\sigma'((x_1,\ldots, x_n,\overline{x}_n,\ldots, \overline{x}_1)) =
(\overline{x}_1,x_2,\ldots, x_n,\overline{x}_n,\ldots, \overline{x}_2,x_1)$.
The proposition implies  $\sigma' S_i = S_{\sigma(i)} \sigma'$,
{}from which the alternative expression $\sigma' =
S_{i_{k+1}} \cdots S_{i_d}S_{\sigma^{ -1}(i_1)} \cdots
S_{\sigma^{-1}(i_k)}$
is available for any $0 \le k \le d$ due to $S_i^2=id$.
We identify $\sigma$ with $\sigma'$, thereby
extend the definition of its domain.
Namely, we also let $\sigma$ act on $B_l$ for any $l$ via
$\sigma =
S_{i_{k+1}} \cdots S_{i_d}S_{\sigma^{-1}(i_1)} \cdots
S_{\sigma^{-1}(i_k)}$.
(We do not exhibit $l$.)
The result is independent of $0 \le k \le d$ and enjoys the
properties $\sigma  \ft_i = \ft_{\sigma(i)}\sigma$ and
$\sigma \et_i = \et_{\sigma(i)}\sigma$.
Proposition \ref{pr:sigma} was known \cite{HKOT} for some $k$.
For the tensor product crystal
$B = B_{l_1} \ot \cd \ot B_{l_N}$, we write
$\sigma_B = \sigma \ot \cd \ot \sigma: B \rightarrow B$, where
$\sigma$ on the right side acts on each component $B_{l_j}$ of the
tensor product according to the above rule
$\sigma =
S_{i_{k+1}} \cdots S_{i_d}S_{\sigma^{-1}(i_1)} \cdots
S_{\sigma^{-1}(i_k)}$.
Obviously one has
$\sigma_B \ft_i = \ft_{\sigma(i)} \sigma_B$ and
$\sigma_B \et_i = \et_{\sigma(i)} \sigma_B$ on $B$, therefore
\begin{equation}\label{eq:Sandsigma}
\sigma_B S_i = S_{\sigma(i)} \sigma_B \quad \text{ on }\;
B = B_{l_1} \ot \cd \ot B_{l_N}.
\end{equation}
The combinatorial $R$ matrix $R: B \ot B' \xrightarrow{\sim} B' \ot B$
satisfies
\begin{equation}\label{eq:new}
R(\sigma_B \ot \sigma_{B'}) = (\sigma_{B'} \ot \sigma_{B})R.
\end{equation}
When acting on crystals,  $\sigma$ without an index shall always act
on a single crystal $B_l$ with some $l$.

For each $a \in \{1,\ldots, n+1\} \,(\hat{\geh}_n = A^{(1)}_{n})$,
$\{1,\ldots, n, \overline{n}, \ldots, \overline{1} \}\,
(\hat{\geh}_n \neq  A^{(1)}_{n})$, we set
\begin{equation}\label{eq:delta}
\delta_l[a] = (x_a = l,\; x_{a'} =0 \, \text{ for }\, a' \neq a) \in B_l
\end{equation}
with the notation in Appendix \ref{app:crystal}.
Using the crystal structure explicitly one easily finds
\begin{align}
&\delta_l[a_k] = S_{i_k}\left(\delta_l[a_{k-1}]\right)
= \et^{\max}_{i_k}\left(\delta_l[a_{k-1}]\right) \quad 1 \le k \le d,
\label{eq:deltaS}\\
&\et^{\max}_{i_{d}} \cd \et^{\max}_{i_{1}} u = \delta_l[a_{d}]\quad
\text{for any } u \in B_l,\label{eq:max}\\
&\varphi_{i_k}(\delta_l[a_{k-1}]) = 0, \nonumber
\end{align}
for any $\hat{\geh}_n$, where $\et^{\max}_ib = \et^{\veps_i(b)}_ib$.
In particular (\ref{eq:deltaS}) implies
$\delta_l[a_d] = \sigma^{-1}(\delta_l[a_0])$ due to
Proposition \ref{pr:sigma}.
Thus it is natural to extend the definition of
$i_k \in I$ and $a_k$
to all $k \in {\mathbb Z}$ by
\begin{equation}\label{eq:extendedai}
i_{k+d} = \sigma^{-1}(i_k),\quad
\delta_l[a_{k+d}] = \sigma^{-1}
\left(\delta_l[a_{k}]\right).
\end{equation}
This way of extending the index $k$ of $a_k$ is independent of $l$, and
(\ref{eq:deltaS}) also persists for all $k \in {\mathbb Z}$.
We have $\{a_k \mid k \in {\mathbb Z} \} =
\{1, \ldots, n+1 \}$ for $\hat{\geh}_n = A^{(1)}_n$
and $\{1, \ldots, n, \overline{n}, \ldots, \overline{1} \}$
for $\hat{\geh}_n \neq A^{(1)}_n$.

In this paper we will concern some asymptotic domain of the
crystal $B_M$ when $M$ gets large.
For $a \in \{a_k \mid k \in {\mathbb Z} \}$
and $M \gg 1$, we introduce the ``domain"
$B_M[a] \subset B_M$ by
\begin{equation}\label{eq:B[]def}
\left.
\begin{array}{ll}
&\hat{\geh}_n = A^{(1)}_{n}:\\
&\;B_M[a] = \{(u_1,\ldots, u_{n+1}) \in B_M \mid u_a \gg u_b 
\text{ for any } b \in \{1,\ldots, n+1\} \setminus \{ a \} \},\\
&\hat{\geh}_n \neq A^{(1)}_{n}:\\
&B_M[a] = 
\left\{ (u_1,\ldots, \overline{u}_1) \in B_M  \Biggm|
\begin{array}{l}
u_a - \overline{u}_a \gg \vert u_b - \overline{u}_b \vert 
\text{ for any } \\
b \in \{1,\ldots, n\}, \; b \neq a, \overline{a}
\end{array}
\right\}.
\end{array}
\right.
\end{equation}
In the rest of the paper we will have assertions under the
conditions like $c \ot u \ot c' \in B \ot B_M[a] \ot B'$
with  $B$ and $B'$ of the form
$B=B_{l_1} \ot \cd \ot B_{l_N}$,
$B'=B_{l'_1} \ot \cd \ot B_{l'_{N'}}$ ($N,N' \ge 0$).
The mathematically invalid ``definition" (\ref{eq:B[]def}) should be
understood that the associated assertions are valid
on condition that $M \gg 1$  and 
the inequalities in (\ref{eq:B[]def}) 
are satisfied. 
Thus for example it amounts to assuming the following:
\begin{align}
&\veps_{i_k} \gg \vphi_{i_k} \; \text{on }
\;B \ot B_M[a_{k-1}]\ot B',\quad
\vphi_{i_k} \gg \veps_{i_k} \; \text{on }
\;B \ot B_M[a_{k}]\ot B',\label{eq:gg} \\
&\veps_{i_{k}}(u \ot x) = \veps_{i_{k}} (u) \;\;\text{ for }
u \ot x \in B_M[a_k] \ot B,\label{eq:eps1} \\
&\vphi_{i_{k+1}}(x \ot u) = \vphi_{i_{k+1}}(u) \;\;\text{ for }
x \ot u \in B \ot B_M[a_k], \label{eq:phi1} \\
&\text{If } \; u \ot x \stackrel{\sim}{\mapsto} y \ot v \; \text{ then }\;
u \ot x \in B_M[a]\ot B \Leftrightarrow y \ot v \in B\ot B_M[a]. \label{eq:iso}
\end{align}
In the above (\ref{eq:gg}) can be checked  by using  the explicit
formula \cite{KKM,HKKOT}
for $\vphi_i$ and $\veps_i$.
(\ref{eq:eps1}) and (\ref{eq:phi1}) follow directly {}from
(\ref{eq:gg}) and (\ref{eq:phi}).
By the weight reason (\ref{eq:iso}) is obvious.
Moreover we may effectively treat as
\begin{equation}
S_{i_k}(B \ot B_M[a_{k-1}]\ot B') =  B \ot B_M[a_k]\ot B',
\label{eq:danger}
\end{equation}
for any $k \in {\mathbb Z}$.

\section{Combinatorial $R$ matrices}\label{sec:CR}

Let $B$ be any crystal of the form
$B = B_{l_1} \ot \cd \ot B_{l_N}$.
Our goal in this section is to prove

\begin{theorem}\label{th:Conj9}
When $M$ is sufficiently large, the combinatorial $R$ matrix
giving the isomorphism $R: B_M \ot B \stackrel{\sim}{\rightarrow} B \ot B_M$
is expressed as
\begin{equation*}
R = (\sigma_B \ot \sigma) P
S_{i_{k+d}}\cd S_{i_{k+2}}S_{i_{k+1}}
\end{equation*}
on $B_M[a_k]\ot B \subset B_M \ot B$ for any $k \in {\mathbb Z}$.
\end{theorem}
\begin{definition}\label{def:Def5}
Let $B$ be any crystal of the form
$B = B_{l_1} \ot \cd \ot B_{l_N}$.
Given $u \ot x \in B_M[a_0] \ot B$ and
$y \ot v \in B \ot B_M[a_0]$ we set
\begin{align*}
\ket{u}{k}\ot \ket{x}{k} &= S_{i_k} \cd S_{i_1}(u\ot x) \in B_M[a_k]\ot B,\\
\bra{y}{k}\ot \bra{v}{k} &= S_{i_k} \cd S_{i_1}(y\ot v)  \in B \ot B_M[a_k],
\end{align*}
for $0 \le k \le d$.
\end{definition}
It should be noted that $\ket{u}{k}$ for example is not defined solely {}from $u$
but only with the other element $x$.

It suffices to show Theorem \ref{th:Conj9} for $k=0$.
To see this, let $u \ot x \in B[a_0]\ot B, \, 1 \le k \le d$ and assume
$u \ot x \stackrel{\sim}{\mapsto} (\sigma_B \ot \sigma)P S_{i_d} \cd S_{i_1}(u \ot x)$
under the isomorphism $B_M \ot B \simeq B \ot B_M$ with $M$ sufficiently large.
Multiplying $S_{i_k} \cd S_{i_1}$ on the both sides and using
(\ref{eq:Sandsigma}), one gets
\begin{equation}\label{eq:mult}
\ket{u}{k} \ot \ket{x}{k} \stackrel{\sim}{\mapsto} (\sigma_B \ot \sigma)P
S_{\sigma^{-1}(i_k)} \cd S_{\sigma^{-1}(i_1)}S_{i_d}\cd S_{i_{k+1}}
(\ket{u}{k} \ot \ket{x}{k}).
\end{equation}
In view of (\ref{eq:extendedai}) and (\ref{eq:danger})
this proves $1 \le k \le d$ case.
To see it for the other $k$,  multiply $(\sigma^m \ot \sigma_B^m)$
on the left side and $(\sigma^m_B \ot \sigma^m)$ on the right side
of (\ref{eq:mult}) for an integer $m$ and use (\ref{eq:Sandsigma}) and
(\ref{eq:new}).

Henceforth we shall concentrate on the  $k=0$ case
of Theorem \ref{th:Conj9} in the rest of this section.
Suppose $B_M[a_0]\ot B \ni u\ot x 
\stackrel{\sim}{\mapsto} y \ot v \in B \ot B_M[a_0]$
under the isomorphism $B_M \ot B \simeq B \ot B_M$
with $M$ sufficiently large.
Then the assertion of Theorem \ref{th:Conj9} with $k=0$ is equivalent to
\begin{align}
v &= \sigma \ket{u}{d}, \label{eq:Conj9.2}\\
y &= \sigma_B \ket{x}{d}.\label{eq:Conj9.1}
\end{align}
We shall prove these relations separately.
In our approach, (\ref{eq:Conj9.2}) can be verified directly for any choice
$B = B_{l_1} \ot \cd \ot B_{l_N}$.
On the other hand, as for (\ref{eq:Conj9.1}), we first deal with
$N=1$ case and derive $N$ general case based on it.

Let us first treat (\ref{eq:Conj9.2}).
\begin{definition}\label{def:Def11}
For any crystal $B$ we introduce
\begin{align*}
t: B &\rightarrow {\mathbb Z}^d_{\ge 0}\\
b &\mapsto (t_1(b), \ldots, t_d(b))\nonumber \\
t_k(b) &= \vphi_{i_k}(\et^{\max}_{i_{k-1}}\cd \et^{\max}_{i_1}b).\nonumber
\end{align*}
\end{definition}
Here $i_k$'s are those in Table \ref{tab:2}.
It is easy to calculate $t$ explicitly for $B = B_l$.
For the element $b = (x_1, \ldots, x_{n+1})$ for $A^{(1)}_n$
($b = (x_1, \ldots, \overline{x}_1)$ for $\hat{\geh}_n \neq A^{(1)}_n$),
the result is summarized in
\begin{lemma}\label{lem:Rem12}
The map $t: B_l \rightarrow {\mathbb Z}^d_{\ge 0}$ has the form:
\begin{align*}
t(b) &= (x_{n},\ldots, x_{1}) \; \text{for }\; A^{(1)}_n,\\
&= (x_n,\ldots, x_3,x_2,{x_1,\overline{x}_1 \brace \overline{x}_1,x_1},
\overline{x}_2, \ldots,\overline{x}_{n-1}) \; \text{for }\; A^{(2)}_{2n-1},\\
&= (x_n,\ldots, x_1,x_0,\overline{x}_1,\ldots,\overline{x}_{n-1}) \; \text{for }\; A^{(2)}_{2n},\\
&= (2x_n+x_0,x_{n-1},\ldots, x_2,{x_1,\overline{x}_1 \brace \overline{x}_1,x_1},
\overline{x}_2, \ldots,\overline{x}_{n-1}) \; \text{for }\; B^{(1)}_{n},\\
&= (x_n,\ldots, x_1,x_0,\overline{x}_1,\ldots,\overline{x}_{n-1}) \; \text{for }\; C^{(1)}_{n},\\
&= (x_n+x_{n-1},x_{n-2},\ldots, x_2,{x_1,\overline{x}_1 \brace \overline{x}_1,x_1},
\overline{x}_2, \ldots,\overline{x}_{n-2},\overline{x}_{n-1}+x_n)
\; \text{for }\; D^{(1)}_{n},\\
&= (2x_n+x_0,x_{n-1},\ldots, x_1,x_{\emptyset},\overline{x}_1,\ldots,\overline{x}_{n-1})
\; \text{for }\; D^{(2)}_{n+1}.
\end{align*}
\end{lemma}
See Appendix \ref{app:crystal} for the notation
$x_0$ in $A^{(2)}_{2n}, C^{(1)}_n$ and $x_{\emptyset}$ in $D^{(2)}_{n+1}$.
The upper and lower choices correspond to those in Table \ref{tab:2}.
{}From Lemma \ref{lem:Rem12} we derive a useful fact.
\begin{proposition}\label{pr:Prop3}
The map $t: B_l \rightarrow {\mathbb Z}^d_{\ge 0}$ is injective for any
$l \in {\mathbb Z}_{\ge 1}$.
\end{proposition}

\begin{lemma}\label{lem:Lem6}
Under Definition \ref{def:Def5} one has
\begin{equation}
\veps_{i_k}(\ket{u}{k}) = \vphi_{i_k}(\ket{u}{k-1}\ot \ket{x}{k-1}),
\quad 1 \le k \le d.
\label{eq:Lem6.1}
\end{equation}
\end{lemma}
{\em Proof}.
Definition \ref{def:Def5} tells
$\ket{u}{k}\ot \ket{x}{k} = S_{i_k}(\ket{u}{k-1}\ot \ket{x}{k-1})$.
Since $\ket{u}{k-1}\ot \ket{x}{k-1} \in B_M[a_{k-1}]\ot B$, this $S_{i_k}$ acts as
\begin{align}
\ket{u}{k}\ot \ket{x}{k} &= \et^q_{i_k}\ket{u}{k-1}\ot \et^{q'}_{i_k}\ket{x}{k-1},
\label{eq:tsukau}\\
q &= \veps_{i_k}(\ket{u}{k-1}) - \vphi_{i_k}(\ket{x}{k-1}) -
\left(\vphi_{i_k}(\ket{u}{k-1}) - \veps_{i_k}(\ket{x}{k-1})\right)_+,\nonumber\\
q' &= \left(\veps_{i_k}(\ket{x}{k-1}) - \vphi_{i_k}(\ket{u}{k-1})\right)_+.\nonumber
\end{align}
Thus we have $\veps_{i_k}(\ket{u}{k}) = \veps_{i_k}(\ket{u}{k-1}) - q
\stackrel{(\ref{eq:phi})}{=}
\vphi_{i_k}(\ket{u}{k-1}\ot \ket{x}{k-1})$.
\qed

For integers $p,q$ depending on $M$ in general, 
we write $p \lesssim q$ to mean $p<q$ or
$p-q = M$-independent for $M \gg 1$.
\begin{lemma}\label{lem:hata}
Let  $u \in B_M[a_0] \subset B_M$ for sufficiently large $M$.
For  $(c_1, \ldots, c_d) \in {\mathbb Z}^d_{\ge 0}$, 
define $u^{\langle k \rangle}\; (0 \le k \le d)$ by
$u^{\langle k \rangle} =
\et^{\veps_{i_k}(u^{\langle k-1 \rangle})-c_k}_{i_k}u^{\langle k-1 \rangle}\;
(1 \le k \le d)$
and $u^{\langle 0 \rangle} = u$.
Suppose $c_j \lesssim t_j(u)$ for all $1 \le j \le d$.
Then the following hold for $1 \le k \le d$:
\begin{align}
&\vphi_{i_k}(u^{\langle k-1 \rangle}) = t_k(u),\label{eq:hata1}\\
&\veps_{i_k}(u^{\langle k \rangle}) = c_k,\label{eq:hata2}\\
&t_k(\sigma u^{\langle d \rangle}) = c_k.\label{eq:hata3}
\end{align}
\end{lemma}
Although $u^{\langle k \rangle}$ depends on $c_j$'s,
the right side of (\ref{eq:hata1}) is independent of them.
Actually (\ref{eq:hata2}) is trivial.
The other relations in
the lemma can be verified with a direct calculation by using the
crystal structure of $B_M$.
Under Definition \ref{def:Def5}, Lemma \ref{lem:hata}  immediately leads to
($1 \le k \le d$)
\begin{align}
&\vphi_{i_k}(\ket{u}{k-1}) = t_k(u),\label{eq:Lem7.1}\\
&\vphi_{i_k}(\bra{v}{k-1}) = t_k(v),\label{eq:Lem7.2}\\
&\veps_{i_k}(\ket{u}{k}) = t_k(\sigma \ket{u}{d}).\label{eq:Lem8}
\end{align}
The right sides of (\ref{eq:Lem7.1}) and (\ref{eq:Lem7.2}) are
independent of $x$ and $y$ in Definition \ref{def:Def5}, respectively.

\vskip0.2cm
{\em Proof of (\ref{eq:Conj9.2})}.
Suppose $B_M[a_0]\ot B \ni u\ot x \stackrel{\sim}{\mapsto} y \ot v \in B \ot B_M[a_0]$
under the isomorphism $B_M \ot B \simeq B \ot B_M$
with $M$ sufficiently large.
We employ the notation in Definition \ref{def:Def5}.
Applying $S_{i_{k-1}}\cd S_{i_1}$ to the both sides of
$u \ot x \stackrel{\sim}{\mapsto} y \ot v$, one gets
$\ket{u}{k-1} \ot \ket{x}{k-1} \stackrel{\sim}{\mapsto} \bra{y}{k-1}\ot \bra{v}{k-1}$,
therefore
$\vphi_{i_k}(\ket{u}{k-1} \ot \ket{x}{k-1})
= \vphi_{i_k}(\bra{y}{k-1}\ot \bra{v}{k-1})$.
But
\begin{align*}
\vphi_{i_k}(\ket{u}{k-1} \ot \ket{x}{k-1}) &
\stackrel{(\ref{eq:Lem6.1})}{=} \veps_{i_k}(\ket{u}{k})
\stackrel{(\ref{eq:Lem8})}{=} t_k(\sigma \ket{u}{d}),\\
\vphi_{i_k}(\bra{y}{k-1}\ot \bra{v}{k-1}) &
\stackrel{(\ref{eq:phi1})}{=}
\vphi_{i_k}(\bra{v}{k-1})
\stackrel{(\ref{eq:Lem7.2})}{=} t_k(v),
\end{align*}
are valid for $1 \le k \le d$, telling that
$t(v) = t(\sigma \ket{u}{d})$.
Thus (\ref{eq:Conj9.2}) follows {}from Proposition \ref{pr:Prop3}.
\qed

Now we proceed to the proof of (\ref{eq:Conj9.1}) with the simple
choice $B = B_l$.

\begin{lemma}\label{lem:Lem14}
Suppose $B_M[a_0]\ot B_l \ni
\delta_M[a_0]\ot z \stackrel{\sim}{\mapsto}
{\tilde z} \ot {\tilde u} \in B_l \ot B_M[a_0]$
under the isomorphism $B_M \ot B_l \simeq B_l \ot B_M$
with $M$ sufficiently large.
Then we have $t({\tilde u}) = t(z)$.
\end{lemma}

{\em Proof}.
Define $\ket{u}{k} \ot \ket{z}{k}$ by Definition \ref{def:Def5}
starting {}from $u \ot z = \ket{u}{0}\ot \ket{z}{0} = \delta_M[a_0] \ot z$.
{}From (\ref{eq:Conj9.2}) we already know that $\tilde{u} = \sigma \ket{u}{d}$.
Thus we have
\begin{align*}
t_k(\tilde{u}) &= t_k(\sigma \ket{u}{d})  \stackrel{(\ref{eq:Lem8})}{=}
\veps_{i_k}(\ket{u}{k}) \stackrel{(\ref{eq:Lem6.1})}{=}
\vphi_{i_k}(\ket{u}{k-1}\ot \ket{z}{k-1}) \\
&\stackrel{(\ref{eq:phi})}{=} \vphi_{i_k}(\ket{z}{k-1}) +
\left( \vphi_{i_k}(\ket{u}{k-1}) - \veps_{i_k}(\ket{z}{k-1}) \right)_+\\
&\stackrel{(\ref{eq:Lem7.1})}{=}
\vphi_{i_k}(\ket{z}{k-1}) +
\left( t_k(u) - \veps_{i_k}(\ket{z}{k-1}) \right)_+.
\end{align*}
Note {}from the explicit forms in Lemma \ref{lem:Rem12} that
$t_k(u) = t_k(\delta_M[a_0]) = 0$, {}from which
$t_k(\tilde{u}) = \vphi_{i_k}(\ket{z}{k-1})$ follows.
In view of $ \vphi_{i_k}(\ket{u}{k-1}) \stackrel{(\ref{eq:Lem7.1})}{=}
t_k(u) = 0$ and (\ref{eq:gg}), it  also follows that
$\ket{z}{k} = \et^{\max}_{i_k} \cd \et^{\max}_{i_1}z$.
Therefore we conclude
$\vphi_{i_k}(\ket{z}{k-1}) = t_k(z)$ {}from Definition \ref{def:Def11}.
\qed

\begin{lemma}\label{lem:Lem15}
Given $y \ot v \in B_l \ot B_M[a_0]$, set
$$
y^{(k)} \ot v^{(k)} = \et^{\max}_{i_k} \cd \et^{\max}_{i_1}(y \ot v)
$$
for $0 \le k \le d$.
Then we have $t(\sigma v^{(d)}) = t(y)$.
\end{lemma}
{\em Proof}.
Since $y^{(k)} = \et^{\max}_{i_k} \cd \et^{\max}_{i_1}y$, we have
$$
v^{(k)} = \et^{\veps_{i_k}(v^{(k-1)})-\vphi_{i_k}(y^{(k-1)})}_{i_k}v^{(k-1)}
= \et^{\veps_{i_k}(v^{(k-1)})-t_k(y)}_{i_k}v^{(k-1)}.
$$
Applying (\ref{eq:hata3}) we obtain $t(\sigma v^{(d)}) = t(y)$.
\qed

\vskip0.2cm
{\em Proof of (\ref{eq:Conj9.1}) for $B = B_l$}.
Suppose $B_M[a_0]\ot B_l \ni u\ot x \stackrel{\sim}{\mapsto} y \ot v \in B_l \ot B_M[a_0]$
under the isomorphism $B_M \ot B_l \simeq B_l \ot B_M$
with $M$ sufficiently large.
We employ the notation in Definition \ref{def:Def5}.
Applying $(\sigma \ot \sigma)\et^{\max}_{i_d} \cd \et^{\max}_{i_1}$
to the both sides and
using (\ref{eq:max}), one finds
$$
B_M[a_0] \ot B_l \ni \delta_M[a_0] \ot \sigma \ket{x}{d}
\stackrel{\sim}{\mapsto}
(\sigma \ot \sigma) \et^{\max}_{i_d} \cd \et^{\max}_{i_1}(y \ot v) \in
B_l \ot B_M[a_0].
$$
Setting
\begin{align*}
\delta_M[a_0] \ot \sigma \ket{x}{d}
& \stackrel{\sim}{\mapsto} {\tilde z}\ot \tilde{u}
\in B_l \ot B_M[a_0],\\
(\sigma \ot \sigma)\et^{\max}_{i_d} \cd \et^{\max}_{i_1}(y \ot v)
&= (\sigma y^{(d)} )\ot ( \sigma v^{(d)}),
\end{align*}
we have $\tilde{z} = \sigma y^{(d)}$ and $\tilde{u} = \sigma v^{(d)}$.
Thus $t(\tilde{u}) = t(\sigma v^{(d)})$.
But we know $t(\tilde{u}) = t(\sigma \ket{x}{d})$ {}from Lemma \ref{lem:Lem14}
and $t(\sigma v^{(d)}) = t(y)$ {}from Lemma \ref{lem:Lem15}.
Therefore $t(\sigma \ket{x}{d}) = t(y)$, compelling
$y = \sigma \ket{x}{d}$ by Proposition \ref{pr:Prop3}.
\qed

To show (\ref{eq:Conj9.1}) for the general choice
$B = B_{l_1} \ot \cd \ot B_{l_N}$, we prepare
\begin{definition}\label{def:Def19}
For $s, s' \in {\mathbb Z}_{\ge 0}, b, b' \in B_l$ and
$i \in I$,
we let  the vertex diagram

\setlength{\unitlength}{1mm}
\begin{picture}(16,13)(-50,0)
\put(3,5){\line(1,0){10}}
\put(8,1){\line(0,1){8}}
\put(-0.1,4.3){$s$}
\put(14,4.3){$s'$}
\put(9,8){$b$}
\put(9,0){$b'$}
\put(6,6){$\scriptstyle{i}$}
\end{picture}
\par\noindent
denote the relations
\begin{equation*}
b' = \et_i^{(\veps_i(b)-s)_+}b, \qquad
s' = \vphi_i(b) + (s-\veps_i(b))_+.
\end{equation*}
\end{definition}
Here $l$ can be any positive integer but we do not exhibit it in the diagram.
The color $i \in I$ is attached to the horizontal line.
The diagram should not be confused with the one
representing the combinatorial $R$ matrix \cite{HHIKTT}.
Given $i$, $(b',s')$ is uniquely fixed {}from $(b,s)$.
Thus for example the diagram

\setlength{\unitlength}{1mm}
\begin{picture}(90,23)(-15,-3)
\put(1,9){$s_0$}
\multiput(6,10)(18,0){2}{\line(1,0){12}}
\multiput(12,5)(18,0){2}{\line(0,1){10}}
\multiput(10,10.9)(18,0){2}{$\scriptstyle{i}$}
\put(20,9){$s_1$}
\put(38,9){$s_2$}
\put(43,10){\line(1,0){5}} \put(49,9){$\cd$}
\put(54,10){\line(1,0){5}} \put(60,9){$s_{N-1}$}
\put(69,10){\line(1,0){12}} \put(75,5){\line(0,1){10}}
\put(83,9){$s_N$}
\put(73,10.9){$\scriptstyle{i}$}

\put(11,16.3){$b_1$}
\put(11,0.7){$b'_1$}

\put(29,16.3){$b_2$}
\put(29,0.7){$b'_2$}

\put(74,16.3){$b_N$}
\put(74,0.7){$b'_N$}

\end{picture}
\par\noindent
is uniquely determined {}from $s_0,\; i$ and $b_1 \ot \cd \ot b_N$.
{}From Definition \ref{def:Def19} it implies
\begin{equation}\label{eq:horizontal}
b'_1 \ot \cd \ot b'_N = \et_i^{\left(\veps_i(b_1 \ot \cd \ot b_N) - s_0\right)_+ }
(b_1 \ot \cd \ot b_N).
\end{equation}

Having established Theorem \ref{th:Conj9} for $B = B_l$, we already know that
\begin{equation*}
B_M[a_0]\ot B_l \ni u \ot b \stackrel{\sim}{\mapsto}
\sigma \ket{b}{d}\ot \sigma \ket{u}{d} \in B_l \ot B_M[a_0]
\end{equation*}
under the isomorphism $B_M \ot B_l \simeq B_l \ot B_M$ for $M$
sufficiently large,
where $\ket{u}{k}\ot \ket{b}{k} = S_{i_k} \cd S_{i_1}(u \ot b)\;
(0 \le k \le d)$.

\begin{proposition}\label{pr:Prop20}
Under the above stated setting, the following diagram holds.

\setlength{\unitlength}{1mm}
\begin{picture}(20,47)(-45,0)
\put(4,10){\line(1,0){12}}
\put(4,25){\line(1,0){12}}
\put(4,35){\line(1,0){12}}
\put(10,5){\line(0,1){10}}
\put(10,20){\line(0,1){20}}
\put(9.5,15.7){$\vdots$}
\put(-4,34){$t_1(u)$} \put(18,34){$t_1(\sigma \ket{u}{d})$}
\put(-4,24){$t_2(u)$} \put(18,24){$t_2(\sigma \ket{u}{d})$}
\put(-4,9){$t_d(u)$} \put(18,9){$t_d(\sigma \ket{u}{d})$}

\put(9,41.2){$b$}
\put(11.3,29){$\ket{b}{1}$}
\put(11.3,20){$\ket{b}{2}$}
\put(11.3,12){$\ket{b}{d-1}$}
\put(9,0.3){$\ket{b}{d}$}

\put(7,36.5){$\scriptstyle{i_1}$}
\put(7,26.5){$\scriptstyle{i_2}$}
\put(7,11.5){$\scriptstyle{i_d}$}
\end{picture}
\end{proposition}

{\em Proof}.
By Definition \ref{def:Def19} we are to check
\begin{align*}
\ket{b}{k} &= \et^{\left(\veps_{i_k}(\ket{b}{k-1}) - t_k(u)\right)_+}_{i_k}
\ket{b}{k-1},\\
t_k(\sigma \ket{u}{d}) &= \vphi_{i_k}(\ket{b}{k-1})
+ \left(t_k(u) - \veps_{i_k}(\ket{b}{k-1})
\right)_+,
\end{align*}
for $1 \le k \le d$.
To show the former, set $x=b$ and
 apply (\ref{eq:Lem7.1}) in $q'$ appearing in (\ref{eq:tsukau}).
The left side of the latter reads
\begin{align*}
t_k(\sigma \ket{u}{d}) &\stackrel{\rm (\ref{eq:Lem8})}{=}
\veps_{i_k}(\ket{u}{k}) \stackrel{\rm (\ref{eq:Lem6.1})}{=}
\vphi_{i_k}(\ket{u}{k-1}\ot \ket{b}{k-1})\\
&\stackrel{(\ref{eq:phi})}{=} \vphi_{i_k}(\ket{b}{k-1}) +
\left( \vphi_{i_k}(\ket{u}{k-1}) - \veps_{i_k}(\ket{b}{k-1})\right)_+,
\end{align*}
which is equal to the right side owing to (\ref{eq:Lem7.1}).
\qed

\vskip0.2cm
{\em Proof of (\ref{eq:Conj9.1}) for $B = B_{l_1} \ot \cd \ot B_{l_N}$}.
Given  any $x = b_1 \ot \cd \ot b_N \in B$ and
$u \in B_M[a_0]$, set $s_{0,k} = t_k(u) \; (1 \le k \le d)$.
Let $\ket{b_j}{k} \in B_{l_j}, s_{j,k} \in {\mathbb Z}_{\ge 0}\;
(1 \le k \le d, 1 \le j \le N)$
be the ones uniquely determined {}from the diagram:

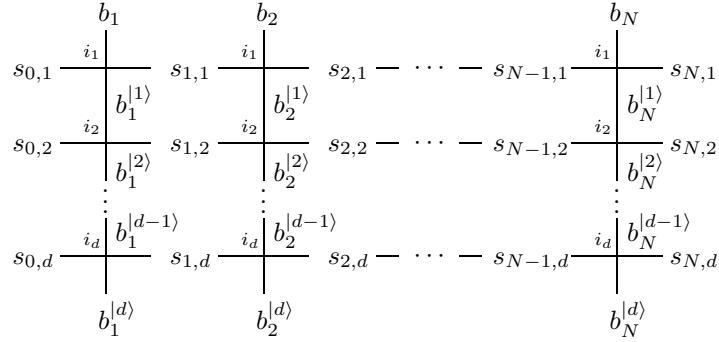
\setlength{\unitlength}{1mm}
\begin{picture}(90,50)(-15,-3)
\put(4,10){\line(1,0){12}}
\put(4,25){\line(1,0){12}}
\put(4,35){\line(1,0){12}}
\put(10,5){\line(0,1){10}}
\put(10,20){\line(0,1){20}}
\put(9.5,15.7){$\vdots$}
\put(-2.4,34){$s_{0,1}$} \put(18.5,34){$s_{1,1}$}
\put(-2.4,24){$s_{0,2}$} \put(18.5,24){$s_{1,2}$}
\put(-2.4,9){$s_{0,d}$}  \put(18.5,9){$s_{1,d}$}

\put(9,41.2){$b_1$}
\put(11.3,29){$\ket{b_1}{1}$}
\put(11.3,20){$\ket{b_1}{2}$}
\put(11.3,12){$\ket{b_1}{d-1}$}
\put(9,0.3){$\ket{b_1}{d}$}

\put(7,36.5){$\scriptstyle{i_1}$}
\put(7,26.5){$\scriptstyle{i_2}$}
\put(7,11.5){$\scriptstyle{i_d}$}

\put(25,10){\line(1,0){12}}
\put(25,25){\line(1,0){12}}
\put(25,35){\line(1,0){12}}

\put(31,5){\line(0,1){10}}
\put(31,20){\line(0,1){20}}
\put(30.5,15.7){$\vdots$}

\put(39.5,34){$s_{2,1}$}
\put(39.5,24){$s_{2,2}$}
\put(39.5,9){$s_{2,d}$}

\put(30,41.2){$b_{2}$}
\put(32.3,29){$\ket{b_{2}}{1}$}
\put(32.3,20){$\ket{b_{2}}{2}$}
\put(32.3,12){$\ket{b_{2}}{d-1}$}
\put(30,0.3){$\ket{b_{2}}{d}$}

\put(28,36.5){$\scriptstyle{i_1}$}
\put(28,26.5){$\scriptstyle{i_2}$}
\put(28,11.5){$\scriptstyle{i_d}$}

\put(46,10){\line(1,0){3}}
\put(46,25){\line(1,0){3}}
\put(46,35){\line(1,0){3}}

\put(51,9.2){$\cdots$}
\put(51,24.2){$\cdots$}
\put(51,34.2){$\cdots$}

\put(57,10){\line(1,0){3}}
\put(57,25){\line(1,0){3}}
\put(57,35){\line(1,0){3}}

\put(61.5,34){$s_{N-1,1}$}
\put(61.5,24){$s_{N-1,2}$}
\put(61.5,9){$s_{N-1,d}$}

\put(72,10){\line(1,0){12}}
\put(72,25){\line(1,0){12}}
\put(72,35){\line(1,0){12}}

\put(78,5){\line(0,1){10}}
\put(78,20){\line(0,1){20}}
\put(77.5,15.7){$\vdots$}

\put(75,36.5){$\scriptstyle{i_1}$}
\put(75,26.5){$\scriptstyle{i_2}$}
\put(75,11.5){$\scriptstyle{i_d}$}

\put(77,41.2){$b_{N}$}
\put(79.3,29){$\ket{b_{N}}{1}$}
\put(79.3,20){$\ket{b_{N}}{2}$}
\put(79.3,12){$\ket{b_{N}}{d-1}$}
\put(77,0.3){$\ket{b_{N}}{d}$}

\put(85,34){$s_{N,1}$}
\put(85,24){$s_{N,2}$}
\put(85,9){$s_{N,d}$}

\end{picture}

\par\noindent
By a repeated use of Proposition \ref{pr:Prop20}, one has
\begin{align*}
u \ot b_1 \ot \cd \ot b_N &\stackrel{\sim}{\mapsto}
\sigma \ket{b_1}{d} \ot \sigma \ket{b_2}{d} \ot \cd
\ot \sigma \ket{b_N}{d} \ot v\\
&= \left(\sigma_B\bigl(\ket{b_1}{d} \ot \cd \ot \ket{b_N}{d}\bigr) \right)
\ot v
\end{align*}
under the isomorphism $B_M \ot B \simeq B \ot B_M$
with sufficiently large $M$.
On the other hand  we introduce
$\ket{u}{k} \ot \ket{x}{k} :=
S_{i_k} \cd S_{i_1}(u \ot x) \in B_M[a_k]\ot B$.
(Although not necessary here, $v \in B_M[a_0]$ is uniquely determined {}from
$t(v) = (s_{N,1}, \ldots, s_{N,d})$ by Proposition \ref{pr:Prop3}, and
we already know that the result coincides with
$v = \sigma \ket{u}{d}$ {}from (\ref{eq:Conj9.2}).)
We are to show
$\ket{x}{d} = \ket{b_1}{d}\ot \cd \ot \ket{b_N}{d}$.
In fact, $\ket{x}{k} = \ket{b_1}{k}\ot \cd \ot \ket{b_N}{k}$
can be proved by induction on $k$ as follows.
(We set $\ket{b_1}{0} \ot \cd \ot \ket{b_N}{0} =
b_1 \ot \cd \ot b_N$.)
It is obvious for $k = 0$.
{}From (\ref{eq:horizontal}) we have
\begin{equation}\label{eq:maruichi}
\ket{b_1}{k}\ot \cd \ot\ket{b_N}{k} =
\et^{(m-s_{0,k})_+}_{i_k}(\ket{b_1}{k-1}\ot \cd \ot \ket{b_N}{k-1}),
\end{equation}
where $m = \veps_{i_k}(\ket{b_1}{k-1}\ot \cd \ot \ket{b_N}{k-1})$.
On the other hand, $\ket{x}{k}$ is determined {}from the recursion relation:
\begin{equation}\label{eq:maruni}
\ket{x}{k} = \et^{\left(\veps_{i_k}(\ket{x}{k-1})-
\vphi_{i_k}(\ket{u}{k-1}) \right)_+}_{i_k} \ket{x}{k-1},
\end{equation}
because of
$\veps_{i_k}(\ket{u}{k-1}) \gg 1$.
Note that
$\vphi_{i_k}(\ket{u}{k-1}) \stackrel{(\ref{eq:Lem7.1})}{=} t_k(u)
= s_{0,k}$.
Thus the two recursion relations (\ref{eq:maruichi}) and (\ref{eq:maruni})
lead to  $\ket{x}{k} = \ket{b_1}{k}\ot \cd \ot \ket{b_N}{k}$
under the induction assumption
$\ket{x}{k-1} = \ket{b_1}{k-1}\ot \cd \ot \ket{b_N}{k-1}$.
\qed

We have finished the proof of (\ref{eq:Conj9.1}) for
any $B = B_{l_1} \ot \cd \ot B_{l_N}$, and thereby the proof of
Theorem \ref{th:Conj9}.

\begin{example}\label{ex:R}
Consider $\hat{\geh}_n = A^{(1)}_3\; (d=3)$.
The data in Table \ref{tab:2} reads
$i_j \equiv a_j \equiv 4-j$ mod 4 with  $i_j \in \{0,1,2,3\}$
and $a_j \in \{1,2,3,4\}$.
To save the space the element
$(x_1,x_2,x_3,x_4) = (3,2,1,0) \in B_6$ is denoted by
$111223$ for example.
Then one has $\sigma(111223) = 112444$ according to Table \ref{tab:1}.

Let us take  $b \otimes c = 111223 \ot 344 \in B_6 \ot B_3$ and seek its
image under the isomorphism
$R: B_6 \ot B_3 \stackrel{\sim}{\rightarrow} B_3 \ot B_6$.
It is known \cite{NaY} that $R(b \ot c) = 223 \ot 111344$.
This can be derived  by taking
$N=1, M=6$ and $k = 3$ in Theorem \ref{th:Conj9}, which reads
$R = (\sigma \ot \sigma) P S_2 S_3 S_0$.
Namely we may regard $b \ot c \in B_M[1]\ot B_3$, where
$M = 6$ and $x_1 = 3$ for $b$ are already sufficiently large so that
\begin{align*}
b \ot c = 111223 \ot 344 &\stackrel{S_0}{\longmapsto} 112234 \ot 344 \\
&\stackrel{S_3}{\longmapsto} 112234 \ot 334 \\
&\stackrel{S_2}{\longmapsto} 112224 \ot 334 \\
&\stackrel{P}{\longmapsto} 334 \ot 112224
\stackrel{\sigma \ot \sigma}{\longmapsto} 223 \ot 111344 = R(b \ot c).
\end{align*}
For comparison, take a smaller example
$R(11223 \ot 344) = 223 \ot 11344 \in B_3 \ot B_5$.
Theorem \ref{th:Conj9} is not applicable in this case
under any choice of $k$
because $S_i(11223 \ot 344) = 11223 \ot 344$ for any $i \in \{0,1,2,3\}$
and $(\sigma \ot \sigma) P(11223 \ot 344) = 233 \ot 11244 \neq 223 \ot 11344$.
\end{example}

\section{Cellular Automata}\label{sec:CA}

The factorization of the combinatorial $R$ matrix
shown in Theorem \ref{th:Conj9} induces the
factorization of the time evolution of the associated cellular automaton.
Consider the isomorphism
\begin{equation*}
B_{M} \ot ( \cdots \ot B_{l_{j}} \ot
B_{l_{j+1}} \ot \cdots )
\xrightarrow{\sim}
( \cdots \ot B_{l_{j}} \ot
B_{l_{j+1}} \ot \cdots ) \ot B_{M}
\end{equation*}
induced by the successive application of combinatorial $R$ matrix
$B_M \ot B_{l_j} \simeq B_{j_j}\ot B_M$.
We impose the boundary condition on $b_j \in B_{l_j}$ as
$b_j = \delta_{l_j}[a_k]$ for $\vert j \vert \gg 1$,
where the choice of $k \in {\mathbb Z}$ is arbitrary.
Assume the following properties:
\begin{align*}
\text{(i)} & \;
\delta_M[a_k] \ot \delta_{l}[a_k] \stackrel{\sim}{\mapsto} \delta_{l}[a_k]\ot \delta_M[a_k]
\quad \text{ for any }\; M, l,\\
\text{(ii)}& \;
u \ot \delta_{l_j}[a_k]\ot \delta_{l_{j+1}}[a_k] \ot \cd \ot \delta_{l_{j+N}}[a_k]
\stackrel{\sim}{\mapsto} \tilde{b}_{j} \ot \cd \ot {\tilde b}_{j+N} \ot \delta_M[a_k]\\
& \text{for any $u \in B_M$ if $N$ is sufficiently large},
\end{align*}
where $\tilde{b}_{l_j} \ot \cd \ot {\tilde b}_{l_{j+N}}$ is some element
in $B_{l_j} \ot \cd \ot B_{l_{j+N}}$.
((i) is indeed valid by the weight reason.)
Then under the isomorphism
$B_M \ot (\cd \ot B_{l_j} \ot B_{l_{j+1}} \ot \cd ) \simeq
(\cd \ot B_{l_j} \ot B_{l_{j+1}} \ot \cd ) \ot B_M$,
\begin{equation}\label{eq:evolution}
\delta_M[a_k] \ot ( \cdots \ot b_{j} \ot
b_{j+1} \ot \cdots )
\stackrel{\sim}{\mapsto}
( \cdots \ot b'_{j} \ot
b'_{j+1} \ot \cdots ) \ot \delta_M[a_k]
\end{equation}
is valid for some $b'_j$'s.
One may regard the system as an automaton which undergoes the time evolution
$p = \cd \ot b_j \ot b_{j+1} \ot \cd \mapsto p' =
\cd \ot b'_j \ot b'_{j+1} \ot \cd$.
When $M$ gets large, the transformation
$p  \mapsto p'$ stabilizes to a certain map,
which we denote \cite{HKT} by $T(p) = p'$.
By taking $B = \cd \ot B_{l_j} \ot B_{l_{j+1}} \ot \cd$ in
Theorem \ref{th:Conj9} and using
(\ref{eq:Sandsigma}), (\ref{eq:extendedai})  we obtain
\begin{corollary}\label{cor:Th21}
Under the boundary condition
$b_j = \delta_{l_j}[a_k] \in B_{l_j}$ for $\vert j \vert \gg 1$,
the time evolution $T$ acts on $B = \cd \ot B_{l_j} \ot B_{l_{j+1}} \ot \cd$ as
\begin{equation*}
T^t = \left\{
\begin{array}{cl}
\sigma^t_B S_{i_{k+td}} \cdots S_{i_{k+2}}S_{i_{k+1}}
& \mbox{ if }\; t \in {\mathbb Z}_{\ge 0},\\
\sigma^t_B S_{i_{k+td+1}}\cd S_{i_{k-1}}S_{i_{k}}
& \mbox{ if }\; t \in {\mathbb Z}_{ < 0}.
\end{array} \right.
\end{equation*}
\end{corollary}
Actually all the  Weyl group operators $S_{i_m}$ for $t > 0$ (resp. $t < 0$)
in the above act as
$S_{i_m} = \et^{\max}_{i_m}$ (resp. $S_{i_m} = \ft^{\max}_{i_m}$)
on $B$ since they  always hit such
states $p \in B$  that  $\veps_{i_m}(p) \gg 1$
(resp. $\vphi_{i_m}(p) \gg 1$).
Corollary \ref{cor:Th21} exhibits a factorization of the time evolution of the
automaton having the  background (vacuum) configuration
$\cd \ot \delta_{l_j}[a_k] \ot \delta_{l_{j+1}}[a_k] \ot \cd $
specified by $k \in {\mathbb Z}$.
When the partial factor
 $S_{i_m}S_{i_{m-1}} \cd S_{i_{k+1}}\, (k+d \ge m \ge k+1)$ in $T$
is applied, the background, hence the boundary condition, changes into
$\cd \ot \delta_{l_j}[a_m] \ot \delta_{l_{j+1}}[a_m] \ot \cd$
according to  (\ref{eq:deltaS}).
A generalization of  $T$
that does not change the background in every intermediate step is
constructed as follows:
\begin{equation}\label{eq:takagi}
\left.
\begin{array}{ll}
{\mathcal T}_m &= \sigma_{B,k,m}
S_{i_m} \cd S_{i_{k+2}} S_{i_{k+1}},\; \; \; \; m \ge k+1,\\
\sigma_{B,k,m} &= \cd \ot \sigma_{k,m} \ot \sigma_{k,m} \ot \cd,\\
\sigma_{k,m} &= S_{i_{k+1}}S_{i_{k+2}} \cd S_{i_m},
\end{array}
\right.
\end{equation}
where $\sigma_{k,m}$ acts on each component $B_{l_j}$ of the tensor product.
We understand ${\mathcal T}_k = id$.
For any $m \ge k$, the operator ${\mathcal T}_m$
retains the background in the original form
$\cd \ot \delta_{l_j}[a_k] \ot \delta_{l_{j+1}}[a_k] \ot \cd $ due to
(\ref{eq:deltaS}).
Moreover {}from Corollary \ref{cor:Th21} we find that the
original evolution under $t$-time step $T^t$ is included in $\{ {\mathcal T}_m \}$
as $T^t = {\mathcal T}_{k+td} \; (t \ge 0)$.
For $A^{(1)}_n$ the evolution of the state $p$ according to
$p = {\mathcal T}_k(p), {\mathcal T}_{k+1}(p), {\mathcal T}_{k+2}(p), \ldots$
agrees with that obtained by the ball-moving algorithm in the
box-ball systems \cite{T,TNS,TTM} under a convention adjustment.
In particular for $A^{(1)}_1$, the evolution rule in terms of the 
Dyck language \cite{TTS} essentially agrees with 
the crystal theoretic signature rule 
in applying $S_i = \et^{\max}_{i}$ or $\ft^{\max}_i\, (i=0,1)$ explained in 
section 1.3 of the reference \cite{KMOU}.

\begin{example}
Consider $\hat{\geh}_n = A^{(2)}_5\; (d=5)$.
To save the space the element
$(x_1,x_2,x_3,\overline{x}_3,\overline{x}_2,\overline{x}_1) = (2,1,2,1,0,1) \in B_7$ is denoted by
$11233\overline{3}\overline{1}$ for example.
Let us take $k=0$, hence the initial background is
$\cd \ot \delta_{l_j}[\,\overline{3}\,] \ot
\delta_{l_{j+1}}[\, \overline{3}\, ]\ot \cd$.
We employ $i_5,\ldots, i_1 = 2,0,1,2,3$ and
 $a_5,\ldots, a_0 = \overline{3}, \overline{2}, 1, 2, 3, \overline{3}$ correspondingly.
Take
$$
p = \cdots \overline{3}\overline{3}\cdot
1\overline{2} \cdot 3 \cdot \overline{3}\overline{1} \cdot 2
\cdot \overline{3}\overline{3} \cdot \overline{3}\overline{3} \cdots
\in \cd \otimes B_2 \ot B_2 \ot B_1 \ot B_2 \ot B_1 \ot B_2 \ot B_2 \ot \cd,
$$
where $\cdot$ stands for $\otimes$ and
$\cdots$ parts on the both ends represent the configuration
identical with the background.
Then the time evolution (\ref{eq:evolution}) is given by
$$
T(p) = \cd \overline{3}\overline{3} \cdot\overline{3}\overline{3} \cdot
 1 \cdot \overline{3}\overline{2} \cdot \overline{1} \cdot 23
\cdot \overline{3}\overline{3} \cd.
$$
According to Corollary \ref{cor:Th21}, the evolution
is decomposed into the following steps.
\begin{align*}
p &= \cdots \overline{3}\overline{3}\cdot
1\overline{2} \cdot 3 \cdot \overline{3}\overline{1} \cdot 2 \cdot \overline{3}\overline{3}
\cdot \overline{3}\overline{3} \cdots\\
S_3(p) &= \cdots 33 \cdot
1\overline{2} \cdot 3 \cdot \overline{3}\overline{1} \cdot 2 \cdot 33 \cdot 33 \cdots\\
S_2S_3(p) &= \cdots 22 \cdot
1\overline{3} \cdot 2 \cdot \overline{3}\overline{1} \cdot 2 \cdot 33 \cdot 22 \cdots\\
S_1S_2S_3(p) &= \cdots 11 \cdot
1\overline{3} \cdot 2 \cdot \overline{3}\overline{2} \cdot 1 \cdot 33 \cdot 11 \cdots\\
S_0S_1S_2S_3(p) &= \cdots \overline{2}\overline{2} \cdot
\overline{3}\overline{2} \cdot \overline{1} \cdot \overline{3}\overline{2} \cdot 1 \cdot 33 \cdot
\overline{2}\overline{2} \cdots\\
S_2S_0S_1S_2S_3(p) &= \cdots \overline{3}\overline{3} \cdot
\overline{3}\overline{3} \cdot \overline{1} \cdot \overline{3}\overline{2} \cdot 1 \cdot 23 \cdot
\overline{3}\overline{3} \cdots\\
\sigma_B S_2S_0S_1S_2S_3(p) &= \cdots \overline{3}\overline{3} \cdot
\overline{3}\overline{3} \cdot 1 \cdot \overline{3}\overline{2} \cdot \overline{1} \cdot 23
\cdot \overline{3}\overline{3} \cdots
= T(p),
\end{align*}
where in the last step we used the fact that $\sigma_B$ interchanges the letters
$1$ and $\overline{1}$ in each component as specified in  Table \ref{tab:1}.
Here the background configuration is changing except the last step.

Alternatively the evolution may also be decomposed
according to (\ref{eq:takagi}) into the
following steps, in which the original background
$\cd \ot \delta_{l_j}[\,\overline{3}\,] \ot
\delta_{l_{j+1}}[\,\overline{3}\,]\ot \cd$ is kept unchanged.
\begin{align*}
p = {\mathcal T}_0(p) &= \cdots \overline{3}\overline{3}\cdot
1\overline{2} \cdot 3 \cdot \overline{3}\overline{1} \cdot 2 \cdot \overline{3}\overline{3}
\cdot \overline{3}\overline{3} \cdots\\
{\mathcal T}_1(p) &= \cdots \overline{3}\overline{3} \cdot
1\overline{2} \cdot \overline{3} \cdot 3\overline{1} \cdot 2 \cdot
\overline{3}\overline{3} \cdot \overline{3}\overline{3} \cdots\\
{\mathcal T}_2(p) &= \cdots \overline{3}\overline{3} \cdot
1\overline{2} \cdot \overline{3} \cdot \overline{2}\overline{1} \cdot \overline{3}
 \cdot 22 \cdot \overline{3}\overline{3} \cdots\\
{\mathcal T}_3(p) &= \cdots \overline{3}\overline{3} \cdot
\overline{3}\overline{2} \cdot 1 \cdot \overline{2}\overline{1} \cdot
\overline{3} \cdot 22 \cdot \overline{3}\overline{3} \cdots\\
{\mathcal T}_4(p) &= \cdots \overline{3}\overline{3} \cdot
\overline{3}\overline{2} \cdot 1 \cdot \overline{3}\overline{2} \cdot \overline{1} \cdot 22 \cdot
\overline{3}\overline{3} \cdots\\
{\mathcal T}_5(p) &= \cdots \overline{3}\overline{3} \cdot
\overline{3}\overline{3} \cdot 1 \cdot \overline{3}\overline{2} \cdot \overline{1} \cdot 23
\cdot \overline{3}\overline{3} \cdots
= T(p).
\end{align*}
\end{example}
Regarding $\overline{3}$ as empty space in a box, one can
interpret the above patterns
as a motion of particles and anti-particles which can form
a neutral (weight 0) bound state.
We hope to report on the explicit algorithm for general
$\hat{\geh}_n$ elsewhere.

\vspace{0.4cm}
\noindent
{\bf Acknowledgements} \hspace{0.1cm}
The authors thank Masato Okado and Yasuhiko Yamada for useful discussions.
One of the authors (A.K.) appreciates the warm hospitality of
the organizers of
{\it The Baxter Revolution in Mathematical Physics},
held at Australian National University, Canberra during February 13--19, 2000,
where a part of this work was presented.

\appendix

\section{Parameterization of $B_l$}\label{app:crystal}

We list the parameterization of the crystal $B_l$.
In  $A^{(1)}_n$ case,
it may be identified with the set of semistandard Young tableaux
of length $l$ one row shape on letters $\{1, \ldots, n+1\}$.
For the other $\hat{\geh}_n$,
$B_l$ may be viewed as a similar set with some additional constraints.
The relevant letters are
$\{1,\ldots, n, \overline{n}, \ldots, \overline{1} \}$ as well as
$0$ and/or $\emptyset$ depending on $\hat{\geh}_n$.
The crystal structure (actions of $\et_i, \ft_i$) can be found in the
article \cite{HKKOT} for $C^{(1)}_n$
and the ones \cite{KKM,KMN2} for the other $\hat{\geh}_n$.

\begin{equation*}
A^{(1)}_n: B_l = \{ (x_{1},\dots ,x_{n+1}) \in
\Z^{n+1}  | x_{i} \geq 0,\,
\sum_{i=1}^{n+1} x_{i} =l \}.\hspace{100pt}
\end{equation*}
%%%%%%%%%%%%%%%%%%%%%%%%%%%%%%%%%%%%%%%%%%%%%%%%%%%%%%%%%%%%
\begin{equation*}
A^{(2)}_{2n-1}: B_l = \{ (x_{1},\dots ,x_{n},\overline{x}_{n},\dots ,\overline{x}_{1}) \in
\Z^{2n}  | x_{i}, \overline{x}_{i} \geq 0,\,
\sum_{i=1}^{n}(x_{i}+\overline{x}_{i}) =l \}.\hspace{100pt}
\end{equation*}
%%%%%%%%%%%%%%%%%%%%%%%%%%%%%%%%%%%%%%%%%%%%%%%%%%%%%%%%
\begin{align*}
A^{(2)}_{2n}:\;
&B_l = \{  (x_1,\dots,x_n,\overline{x}_n,\dots,\overline{x}_1) \in
\Z^{2n} \;| \; x_i,\overline{x}_i \ge 0,\, \sum_{i=1}^n (x_i + \overline{x}_i)
\le l \}.\hspace{100pt}\\
&\text{We set }x_0 = l - \sum_{i=1}^n (x_i + \overline{x}_i).
\end{align*}
%%%%%%%%%%%%%%%%%%%%%%%%%%%%%%%%%%%%%%%%%%%%%%%%%%%%%%%%%%%%%%%%%%%%%%%%%
\begin{equation*}
B^{(1)}_{n}:
B_l = \left\{ (x_{1},\dots ,x_{n},x_{0},\overline{x}_{n},\dots ,\overline{x}_{1}) \in
\Z^{2n} \times\{0,1\}
\left|
\begin{array}{c} x_{0} = 0\; \mbox{or} \; 1, x_{i}, \overline{x}_{i} \geq 0, \\
x_0 + \sum_{i=1}^{n}(x_{i}+\overline{x}_{i}) =l
\end{array}\right.
\right\}.\hspace{100pt}
\end{equation*}
%%%%%%%%%%%%%%%%%%%%%%%%%%%%%%%%%%%%%%%%%%%%%%%%%%%%%%%%%%%%%%%%%%%%%%%%
\begin{align*}
C^{(1)}_{n}:\;
&B_l  =  \left\{  (x_1,\dots,x_n,\overline{x}_n,\dots,\overline{x}_1) \in
\Z^{2n} \Biggm|
\begin{array}{l}
x_i,\overline{x}_i \ge 0, \\
l \ge \sum_{i=1}^n
(x_i + \overline{x}_i) \in l - 2 \Z
\end{array}
\right\}.\hspace{100pt}\\
&\text{We set }x_0 = (l-\sum_{i=1}^n(x_i + \overline{x}_i))/2.
\end{align*}
%%%%%%%%%%%%%%%%%%%%%%%%%%%%%%%%%%%%%%%%%%%%%%%%%%%%%%%%%%%%%%%%%%
\begin{equation*}
D^{(1)}_{n}:
B_l = \left\{ (x_{1},\dots ,x_{n},\overline{x}_{n},\dots ,\overline{x}_{1}) \in
\Z^{2n}  \left|
\begin{array}{c}
x_{n} = 0\; \mbox{or} \;\overline{x}_{n} =0,\; x_{i}, \overline{x}_{i} \geq 0,\\
\sum_{i=1}^{n}(x_{i}+\overline{x}_{i}) =l
\end{array}\right.
\right\}.\hspace{100pt}
\end{equation*}
%%%%%%%%%%%%%%%%%%%%%%%%%%%%%%%%%%%%%%%%%%%%%%%%%%%%%%%%%%%%%%%%%%%%%%
\begin{align*}
D^{(2)}_{n+1}:\;
&B_l = \left\{  (x_1,\dots,x_n,x_0,\overline{x}_n,\dots,\overline{x}_1) \in
\Z^{2n} \times \{0,1\} \Biggm|
\begin{array}{l}
x_0=\mbox{$0$ or $1$},x_i,\overline{x}_i \ge 0, \\
x_0 + \sum_{i=1}^n (x_i + \overline{x}_i) \le l
\end{array}
\right\}.\hspace{100pt}\\
&\text{We set }x_{\emptyset} = l-x_0 - \sum_{i=1}^n (x_i + \overline{x}_i).
\end{align*}

%%%%%%%%%%%%%%%%%%%%%%%%%%%%%%%%%%%%%%%%%%%%

\end{document}